\documentclass[11pt]{MBE}

\usepackage{epsf, subfigure, verbatim, epsfig}

\usepackage{colordvi,fancybox,latexsym}
\usepackage{graphics}
\usepackage{graphicx}
\usepackage{epic}
\usepackage{color}

\def\currentvolume{?}
\def\currentissue{?}
\def\currentyear{200?}
\def\currentmonth{?}
\def\ppages{\pageref{title}--\pageref{end}}

\newtheorem{theorem}{Theorem}
\newtheorem{remark}[theorem]{Remark}
\newtheorem{definition}[theorem]{Definition}

\title[Stability in the Perturbed Chemostat]
      {On the Stability of Periodic Solutions\\ in the Perturbed
Chemostat}
      \label{title}
\author[F. Mazenc, M. Malisoff and  P. De Leenheer]{}


 \email{Frederic.Mazenc@ensam.inra.fr}
 \email{malisoff@lsu.edu}
 \email{deleenhe@math.ufl.edu}
\urladdr{http://www.math.lsu.edu/$\sim$malisoff/}
\urladdr{http://www.math.ufl.edu/$\sim$deleenhe/}

\begin{document}

\maketitle

\centerline{\scshape   Fr\'ed\'eric Mazenc}
 \smallskip

  {\footnotesize \centerline{Projet MERE
INRIA-INRA}\centerline{UMR Analyse des Syst\`{e}mes et Biom\'{e}trie
INRA}\centerline{ 2, pl. Viala,
        34060 Montpellier, France}}\medskip\smallskip

\centerline{\scshape Michael Malisoff}
 \smallskip

  {\footnotesize \centerline{Department of Mathematics}\centerline{Louisiana
  State University}
  \centerline{Baton Rouge, LA 70803-4918 } }
  \medskip\smallskip

\centerline{\scshape Patrick De Leenheer}
 \smallskip

  {\footnotesize \centerline{Department of Mathematics}\centerline{University of Florida}\centerline{411
Little Hall, PO Box 118105}\centerline{Gainesville, FL 32611--8105}}

\medskip\smallskip

\centerline{(Communicated by ????)}

 \medskip


\begin{abstract}
We study the chemostat model for one species competing for one
nutrient using a Lyapunov-type analysis. We design the dilution rate
function so that all solutions of the chemostat converge to a
prescribed periodic solution.
 In terms of chemostat biology, this means that no matter what
positive initial levels for the species concentration and nutrient
are selected, the long term species concentration and substrate
levels closely approximate a prescribed oscillatory behavior.  This
is significant because it reproduces the realistic ecological
situation where the species and substrate concentrations oscillate.
We show that the stability is maintained when the model is augmented
by additional species that are being driven to extinction. We also
give an input-to-state stability result for the chemostat tracking
equations for cases where  there are small perturbations acting on
the dilution rate and initial concentration. This means that the
long term  species concentration and substrate behavior enjoys a
highly desirable robustness property, since it continues to
approximate the prescribed oscillation up to a small error when
there are small unexpected changes in the dilution rate function.
\end{abstract}

\smallskip\medskip

\section{Introduction}
The chemostat model provides the foundation for much of current
research in bio-engineering, ecology, and  population biology
\cite{DLLS03,DLS03, EPR01, GR05, SW95}. In the engineering
literature, the chemostat is known as the continuously stirred tank
reactor. It has been used for modeling the dynamics of interacting
organisms in waste water treatment plants, lakes and oceans. In its
basic setting, it describes the dynamics of species competing for
one or more limiting nutrients. If there are $n$ species with
concentrations $x_i$ for $i=1,\dots, n$ and just one limiting
nutrient with concentration $S$ and dilution rate $D>0$, then the
model takes the form

\begin{equation} \label{model-full}
\left\{\begin{array}{rcl}
{\dot S}&=&D(S_{in}-S)-\displaystyle\sum_{i=1}^n\mu_i(S)x_i/\gamma_i\\
{\dot x_i}&=&x_i(\mu_i(S)-D),\;\; i=1,\dots, n
\end{array}\right.
\end{equation}

\begin{figure}[h]
\centerline{\scalebox{.77}{\input{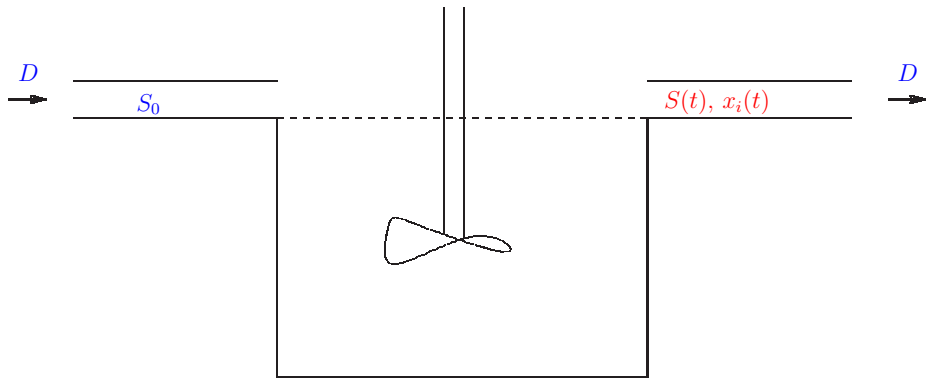}}}
\caption{Chemostat} \label{chem}
\end{figure}

\smallskip\smallskip

\noindent where $\mu_i$ denotes the per capita growth rate of
species $i$ and $\dot p$ is the time derivative of any variable $p$.
(In much of the paper, we simplify our notation by omitting the
arguments of the functions. For instance, when no confusion can
arise from the context, we denote $S(t)$  simply by $S$.) The
functions $\mu_i$ depend only on the nutrient concentration, and are
zero at zero, continuously differentiable and strictly increasing,
although non-monotone functions have been the subject of research as
well. The conversion of nutrient into new biomass for each species
$i$ happens with a certain yield $\gamma_i\in(0,1)$ and the natural
control variables in this model are the input nutrient concentration
$S_{in}$ and the dilution rate $D$. The latter variable is defined
as the ratio of the volumetric flow rate $F$ (with units of volume
over time) and the reactor volume $V_r$ which is kept constant.
Therefore it is proportional to the speed of the pump that supplies
the reactor with fresh medium containing the nutrient. The equations
(\ref{model-full}) are then straightforwardly obtained from writing
the mass-balance equations for the total amounts of the nutrient and
each of the species, assuming the reactor content is well-mixed. The
full model (\ref{model-full}) is illustrated in Figure \ref{chem}.

In the present work, we consider the case where there is  just one
species with concentration $x$, in which case the equations
(\ref{model-full}) take the form
\begin{equation} \label{model-reduced}
\left\{\renewcommand{\arraystretch}{1.25}\begin{array}{rcl}
{\dot S}&=&D(S_{in}-S)-\mu(S)x/\gamma\\
{\dot x}&=&x(\mu(S)-D)
\end{array}\right.
\end{equation}
\noindent(but see Theorem \ref{iss-track} below for results on
chemostats with disturbances, and Section \ref{several} for models
involving several species). We assume  $S_{in}$ is a given positive
constant, while the per capita growth rate $\mu$ is a Monod function
(which is also known as a Michaelis-Menten function) taking the form
\begin{equation}\label{muchoice}
\mu(S)=\dfrac{mS}{a+S},
\end{equation}
\noindent for certain positive constants $m$ and $a$ that we specify
later.  The dilution rate is an appropriate continuous positive
periodic function  we also specify below.  Since $\dot S\ge 0$ when
$S> 0$ is near zero, one  can readily check that
(\ref{model-reduced}) leaves the domain of interest
\[\mathcal{X}:=(0,\infty)\times(0,\infty)\] positively invariant (i.e.,
trajectories for (\ref{model-reduced}) starting in $\mathcal{X}$
remain in $\mathcal{X}$ for all future times);  see Theorem
\ref{iss-track} for a more general invariance result for perturbed
chemostats.

Since we are taking $S_{in}$ to be a fixed positive constant, we
rescale the variables to reduce the number of parameters. Using the
change of variables
\[
{\bar S}=\dfrac{S}{S_{in}}, \; \; {\bar x}=\dfrac{x}{S_{in}\gamma}, \; \;
{\bar \mu}({\bar S})=\mu(S_{in}{\bar S})
\]
\noindent and dropping bars, we eliminate  $S_{in}$ and $\gamma$ and
so obtain the new dynamics
\begin{equation} \label{model}
\left\{\renewcommand{\arraystretch}{1.35}\begin{array}{rcl} {\dot
S}&=&D(1-S)-\mu(S)x\\ {\dot x}&=&x(\mu(S)-D)\end{array}\right.
\end{equation}
again evolving on the state space $\mathcal{X}=(0,\infty)^2$.
\textcolor{black}{Motivated by the realistic ecological situation
where the species concentrations are oscillating, we solve the
following biological problem:\begin{itemize}
\item[]\textcolor{blue}{{\bf Biological Problem B1:}}
For a prescribed oscillatory behavior for the species concentration
and substrate level given by a time-periodic pair $(S_r(t),
x_r(t))$, design a dilution rate function $D(t)$  such that if this
choice of $D(t)$ is used in the chemostat (\ref{model}), then all
solution pairs $(S(t),x(t))$ for the substrate levels and
corresponding
 species levels obtained from solving (\ref{model}) (i.e. for all possible initial
values) closely approximate $(S_r(t),x_r(t))$ for large times $t$.
\end{itemize}See also Problem (SP) in Section \ref{define} below for a precise mathematical statement of
the preceding problem. In the language of control theory, solving
Biological Problem B1  means we will prove the stability of a
suitable periodic reference signal for the species concentration
$t\mapsto x_r(t)$ in (\ref{model}) using an appropriate
time-periodic dilution rate $D(t)$;  see \cite{K02} for the
fundamental ideas from control theory we need in the sequel.}

\textcolor{black}{Since $D(t$) is proportional to the speed of the
pump which supplies the chemostat with medium containing nutrient,
implementation of the prescribed oscillatory behavior requires that
we control the pump in a very precise way. In practice this control
process is prone to errors, and the actual pump speed will be
subject to small fluctuations which we  will model by replacing
$D(t)$ by $D(t)+u_1(t)$ in the chemostat equations, where $u_1(t)$
models the error. It is therefore of interest to study the effect of
these small fluctuations on the periodic behavior. Preferably this
effect will be small, and the resulting behavior is not too
different from the prescribed periodic behavior. We will show that
this is indeed the case by actually quantifying how small the
deviations are, relying on the well-known control-theoretic notion
of Input-to-State Stability or ISS;  see \cite{S00, S06} and Remark
\ref{aboutISS} for details about ISS.   Summarizing this a bit more
formally, we will solve the following biological problem (which we
state in a more precise mathematical way in Section \ref{thm}):}
\textcolor{black}{\begin{itemize}
\item[]\textcolor{blue}{{\bf Biological Problem B2:}}  For the prescribed oscillatory
 behavior $(S_r(t),x_r(t))$ and dilution rate
$D(t)$ obtained in Biological Problem B1,  quantify how the
substrate and species levels $(S(t),x(t))$ in the chemostat model
(\ref{model}) are affected by unexpected changes in the dilution
rate, and  show that the convergence of $(S(t),x(t))$ to the
oscillatory behavior $(S_r(t),x_r(t))$ is robust to small changes in
$D(t)$.\end{itemize}   Our solution to Biological Problem B2 will be
a special case of our more general input-to-state stability result
for the chemostat tracking equations, assuming the dilution rate and
initial concentration are both perturbed by  {}noise terms of small
enough magnitude}.

In the next section, we briefly review the literature focusing on
what makes our approach different.  In Section \ref{track}, we fix
the reference signal we wish to track.  In Section \ref{define},
 we precisely formulate the  definitions and the stability
problem we are solving.   We state our main stability theorem in
Section \ref{thm} and we discuss the significance of our theorem in
Section \ref{discuss}. We prove our
 stability result in Section
\ref{main}.  In Section \ref{several}, we show that the stability is
maintained when there are additional species that are being driven
to extinction.  We validate our results  in Section
\ref{simulations} using a numerical example. We conclude in Section
\ref{concl} by summarizing our findings.

\section{Review of the Literature and Comparison with Our Results}
\label{review} The behavior of the system $(\ref{model-full})$ is
well understood when $S_{in}$ and $D$ are positive constants, as
well as cases where $n=2$ and either of these control variables is
held fixed while the other is periodically time-varying. See
\cite{HS83, S81} for periodic variation of $S_{in}$ and \cite{BHW85}
for periodic variation of $D$ and the general reference  \cite{SW95}
on chemostats. When both $S_{in}$ and $D$ are constants, the
so-called ``competitive exclusion principle'' holds, meaning that at
most one species survives. Mathematically this translates into the
statement that system $(\ref{model-full})$ has a steady state with
at most one nonzero species concentration, which attracts almost all
solutions; see \cite{SW95}. This result has triggered much research
to explain the discrepancy between the (theoretical) competitive
exclusion principle and the observation that in real ecological
systems, many species coexist.

The results on the periodically-varying chemostat mentioned above
should be seen as attempts to explain this paradox. They involve
chemostats with $n=2$ species, and their purpose is to show that an
appropriate periodic forcing for either $S_{in}(t)$ or $D(t)$ can
make the species coexist, usually in the form of a (positive)
periodic solution. Few results on coexistence of $n>2$ species are
available. An exception is \cite{RR90}, where a periodic function
$S_{in}(t)$ is designed (with $D$ kept fixed) so that the resulting
system has a (positive) periodic solution with an arbitrary number
of coexisting periodically varying species. The stability properties
of this solution are not known.

More recent work has explored the use of state-dependent but time
invariant feedback control of the dilution rate $D$ to generate
coexistence; see \cite{DLLS03,DLS03} for monotone growth rate
functions in the $n=2$ species case, and \cite{DLP05} for the $n=3$
species case. The paper \cite{GR05} considers feedback control when
the growth rate functions are non-monotone. In  \cite{GMR05},
\cite{LMR1},  and \cite{MLR}, coexistence is proved for models
taking into account intra-specific competition. In these models, the
usual growth functions $\mu_i(S)$ are replaced by functions
$\mu_i(S,x_i)$ which are decreasing with respect to the variable
$x_i$. All the results discussed so far apply to a more general
model than $(\ref{model})$  involving $n>1$ species. This is because
the main purpose of these papers is to investigate environmental
conditions under which the competitive exclusion principle fails and
several species can coexist.

Here we will not consider any coexistence problems. Our main
objective is to provide a proof of stability of a periodic solution
based on a Lyapunov-type analysis and to investigate the robustness
properties of the periodic solution with respect to perturbations.
As an illustration we show that the stability of the periodic
solution is robust with respect to additional species that are being
driven to extinction, or to small disturbances on the initial
nutrient concentration or dilution rate. These features set our work
apart from the known results on periodically forced chemostat models
which do not rely on the construction of a Lyapunov function.
Proving stability in the chemostat usually relies on reduction and
monotonicity arguments, and not so often on Lyapunov functions (but
see for instance Theorem $4.1$ in \cite{SW95} which uses a Lyapunov
function introduced in \cite{H78} and more recently \cite{GMR05}).

Finally we point out that closely related to our results is
\cite{EPR01} where a single-species chemostat with a continuous and
bounded (but otherwise arbitrary) function $S_{in}(t)$ and constant
dilution rate is investigated; there it is shown that two positive
solutions converge to each other. However, the proof is not based on
a Lyapunov function.  The advantage of having a Lyapunov function is
that it can be used to {\em quantify} the effect of additional noise
terms on the stability of the unperturbed dynamics. In fact, to our
knowledge, our work provides the first input-to-state stability
analysis of chemostats whose dilution rates and initial
concentrations are perturbed by small noise; see Remark
\ref{aboutISS} for a discussion on the importance of input-to-state
stability in control theory and engineering applications.

\section{Choosing a Reference Trajectory} \label{track} 
\begin{figure}[h]
\begin{center}
\scalebox{.6}{\includegraphics{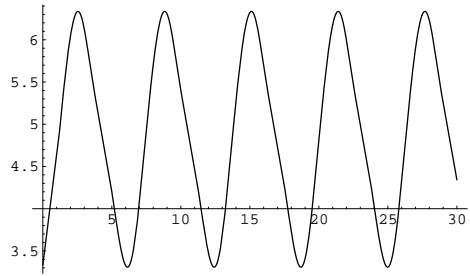}}
\end{center}
\caption{Dilution Rate  $D(t)$ for the Chemostat from (\ref{chli})
\label{refff}}
\end{figure}

We first
choose   the dilution rate $D=D(t)$ that will give  a  reference
trajectory $(S_r(t),x_r(t))$ for (\ref{model}) which we show to be
stable. We assume a growth rate with constants $m,a>0$ as
\textcolor{black}{follows, in which the constants $a$ and $m$ and
the variable $S$ are all dimensionless by the change of coordinates
used to obtain the normalized equations (\ref{model}) so the units
do not matter:}
\begin{equation}
\label{mbound} \mu(S)=\frac{mS}{a+S},\; \; {\rm where}\; \;
 m>4a+1\; .
\end{equation}
\textcolor{black}{For the sake of computational simplicity, we
choose a sinusoidal reference trajectory but the extension to more
general reference trajectories can be handled by similar methods;
see Remark \ref{means} for details.  Simple calculations show that
(\ref{model}) admits the trajectory}
\begin{equation}
\label{reftraj} (S_r(t),x_r(t)) =
\left(\frac{1}{2}-\frac{1}{4}\cos(t),\frac{1}{2}+\frac{1}{4}\cos(t)\right)
\end{equation}
which we refer to as a {\em reference trajectory}  when we choose
\begin{equation}
\label{chli} D(t)\; \; =\; \; -\frac{\dot
x_r(t)}{x_r(t)}+\mu(1-x_r(t))\; \; =\; \;
\frac{\sin(t)}{2+\cos(t)}+\frac{m(2-\cos(t))}{4a+2-\cos(t)}.
\end{equation}
 Condition (\ref{mbound}) then provides constants
$\bar D, D_o> 0$  such that \[D_o\; \le\;  D(t)\; \le\;  \bar D\]
for all $t \ge 0$:
\begin{equation}
\label{aer1} \bar D = 1 + \frac{3m}{4a + 3} \; \; {\rm and}\;  \;
D_o = \frac{m}{4a + 1}-1\; .
\end{equation}
See Figure    \ref{refff} for the graph of  $D(t)$ for $m=10$ and
$a=\frac{1}{2}$.

\section{Definitions and Statement of Stability Problem}
\label{define}

We wish to solve the following stability problem
\textcolor{black}{which is merely a restatement of Biological
Problem B1 above in precise control theoretic terms}:

\begin{itemize}
\item[]\begin{itemize}
\item[(SP) \ ]
Given any trajectory $(S,x): [0,\infty) \to (0,\infty)^2$ for
(\ref{model}) corresponding to the dilution rate  $D(t)$ from
(\ref{chli}) and $\mu$ as in (\ref{mbound}) (i.e. for any initial
value for $(S,x)$), show that the corresponding deviation $(\tilde
S(t), \tilde x(t)):=(S(t) - S_r(t), x(t)-x_r(t))$ of $(S,x)$ from
the reference trajectory (\ref{reftraj}) asymptotically approaches
$(0,0)$ as $t\to +\infty$.
\end{itemize}\end{itemize}

 We will
solve (SP) by proving a far more general tracking result for a
single species chemostat acted on by a disturbance vector
$u=(u_1,u_2):[0,\infty)\to \mathbb{R}^2$ as follows:
\begin{equation}
\label{mod1} \left\{
\begin{array}{rcl}
\dot{S}(t) & = & [D(t) +
u_1(t)](1+u_2(t)-S(t))-\mu(S(t))x(t)\\[.5em]
\dot x(t)&=& x(t)[\mu(S(t))-D(t)-u_1(t)]
\end{array}
\right..
\end{equation}
We will quantify the extent to which the reference trajectory
(\ref{reftraj}) tracks the trajectories of (\ref{mod1})
\textcolor{black}{which will solve Biological Problem B2 from the
introduction}. To this end, we need to introduce a priori bounds on
$u_1$ and $u_2$; see Remark \ref{bound-u}. Our main theoretical tool
will be the input-to-state stability (ISS) property \cite{S89} which
is one of the central paradigms  of current research in nonlinear
stability analysis;  see Remark \ref{aboutISS}. The relevant
definitions are as follows.

We let ${\mathcal K}_\infty$ denote the set of all continuous
functions $\gamma:[0,\infty)\to[0,\infty)$ for which (i)
$\gamma(0)=0$ and (ii) $\gamma$ is strictly increasing and
unbounded. We let $\mathcal{KL}$ denote the class of all continuous
functions $\beta:[0,\infty)\times [0,\infty)\to[0,\infty)$ for which
(I)
 $\beta(\cdot, t)\in {\mathcal K}_\infty$ for each $t\ge 0$,
 (II)
  $\beta(s,\cdot)$ is non-increasing for each $s\ge 0$, and
(III)  $\beta(s,t)\to 0$ as $t\to +\infty$ for each $s\ge 0$.
Consider a general control-affine dynamic
\begin{equation}
\label{gen} \dot y=F(y,t)+G(y,t)u,\; \; y\in \mathcal{O},\; u\in
\mathbf{U}
\end{equation}
evolving on a given  subset $\mathcal{O}\subseteq\mathbb{R}^n$ where
$\mathbf U$ is a given subset of Euclidean space. (Later we
specialize to dynamics for the chemostat.)  For each $t_o\ge 0$ and
$y_o\in \mathcal{O}$, let $y(t; t_o, y_o, \alpha)$ denote the
solution of (\ref{gen}) satisfying $y(t_o)=y_o$ for a given control
function $\alpha\in \mathcal{U}:=\{{\rm measurable\ essentially\
bounded\ }\alpha:[0,\infty)\to {\mathbf U}\}$; i.e. the solution of
the initial value problem \[\dot y(t)=F(y(t),t)+G(y(t),t)\alpha(t)
\; {\rm a.e.}\; t\ge t_o\, ,\; y(t_o)=y_o\; .\]
  We always assume that such solutions are uniquely defined on
all of $[t_o,\infty)$ (i.e., (\ref{gen}) is forward complete and
$\mathcal{O}$ is positively invariant for this system) and that
there exists $\Theta\in \mathcal{K}_\infty$ such that
$|F(y,t)|+|G(y,t)|\le \Theta(|y|)$ everywhere, where $|\cdot |$ is
the usual Euclidean norm.  For example,
\[[t_o,\infty)\ni t\mapsto (S(t; t_o,(S_o,x_o),\alpha),x(t; t_o,
(S_o,x_o),\alpha))\] is the solution of (\ref{mod1}) for the
disturbance $u=(u_1, u_2)=\alpha(t)$ satisfying the initial
condition $(S(t_o),x(t_o))=(S_o,x_o)$.

\begin{definition}\label{issdef}
We call (\ref{gen})  {\em input-to-state stable (ISS)} provided
there exist $\beta\in \mathcal{KL}$ and $\gamma\in
\mathcal{K}_\infty$ such that
\begin{equation}
\label{ISSest} |y(t; t_o, y_o, \alpha)|\; \; \le\; \;  \beta(|y_o|,
t-t_o)+\gamma(|\alpha|_\infty)
\end{equation}
for all $t\ge t_o$, $t_o\ge 0$, $y_o\in \mathcal{O}$, and $\alpha\in
\mathcal{U}$.
\end{definition}
Here $|\alpha|_\infty$ denotes the essential supremum of $\alpha\in
\mathcal{U}$.  By causality, the ISS condition (\ref{ISSest}) is
unchanged if $|\alpha|_\infty$ is replaced by the essential supremum
$|\alpha|_{[t_o,t]}$ of $\alpha$ restricted to $[t_o,t]$.  In
particular, (\ref{ISSest}) says $y(t; t_o, y_o, \mathcal{Z})\to 0$
as $t\to +\infty$ for all initial values $y_o$ and initial times
$t_o$, where $\mathcal{Z}$ is the zero disturbance $\alpha(t)\equiv
0$.

\begin{remark}
\label{aboutISS}  The theory of ISS systems originated in
\cite{S89}. ISS theory provides the foundation for much current
research in robustness analysis and controller design for nonlinear
systems, and has also been used extensively in engineering and other
applications \cite{A04, AISW04, C05, MRS04, S89, S06}. The ISS
approach can be viewed as a unification of the operator approach of
Zames (e.g. \cite{Z66a, Z66b}) and the Lyapunov state space
approach. The operator approach involves studying the mapping $(t_o,
y_o, \alpha)\mapsto y(\cdot ; t_o, y_o,\alpha)$ of initial data and
control functions into appropriate spaces of trajectories, and  it
has the advantages that it allows the use of Hilbert or Banach space
techniques to generalize many properties of linear systems to
nonlinear dynamics. By contrast, the state space approach is well
suited to nonlinear dynamics and lends itself to the use of
topological or geometric ideas. The ISS framework has the advantages
of both of these approaches including an equivalent characterization
in terms of the existence of suitable Lyapunov-like functions; see
Remark \ref{pba} below. For a comprehensive survey on many recent
advances in ISS theory including its extension to systems with
outputs, see \cite{S06}.

\end{remark}

To specify the bound $\bar u$ on our disturbances $u=(u_1,u_2)$, we
use the following constants whose formulas will be justified by the
proof of our main stability result:
\begin{equation}
\label{ck}
\begin{array}{l}
c=8\left(\frac{1}{2}+\frac{\bar D}{D_o}\right)^2,\; \; \; \kappa=4+
\max\left\{\frac{112m}{4a+1}\, ,\,
\frac{16m(4a+3)(a+2)}{a(4a+1)^2}\right\}\;
,\\
C_1=\min\left\{1,\frac{\kappa}{200},\frac{ma}{2(4a+3)(a+2)}\right\}
\end{array}
\end{equation}

\section{Statement of Theorem}
\label{thm}

{}From now on, we assume the disturbance vector $u=(u_1,u_2)$ in
(\ref{mod1}) takes all of its values in a fixed square control set
of the form

\begin{equation}
\label{constraint} \begin{array}{l} \mathbf{U}\; :=\; \{(u_1,u_2)\in
\mathbb{R}^2: \, |u_1|\le \bar u, \, |u_2|\le \bar u\} \; \; \text{
where}\\[.5em]
 0 < \bar u\; <\; \min\left\{\dfrac{C_1}{\sqrt{ 8(1 +
2c\kappa C_1)}}, \dfrac{D_o}{2}\right\}
\end{array}
\end{equation}
where $c$,  $\kappa$, and $C_1$ are in (\ref{ck}) (but see Remark
\ref{bound-u} for related results under less stringent conditions on
the disturbance values). We will prove the following robustness
result:
\begin{theorem}
\label{iss-track}  Choose $D(t)$, $\mu$, and $(S_r,x_r)$ as in
(\ref{mbound})-(\ref{chli}).  Then the corresponding solutions of
(\ref{mod1}) satisfy
\begin{equation}
\label{c1ru1}
\renewcommand{\arraystretch}{1.25}\begin{array}{l}\left[\,  S_o> 0\; \;
\&\; \; x_o>0\; \; \&\; \;  t\ge t_o\ge 0\; \; \&\; \; \alpha\in
\mathcal{U}\, \right]\\[.25em]\; \; \; \Rightarrow\; \; \;  \left[\; S(t; t_0,
(S_0,x_0), \alpha)> 0 \; \; \& \; \; x(t; t_0, (S_0,x_0),
\alpha)>0\; \right]\; \; .\end{array}
\end{equation}
Moreover, there exist $\beta\in \mathcal{KL}$ and $\gamma\in
\mathcal{K}_\infty$ such that the corresponding transformed error
vector \[\renewcommand{\arraystretch}{1.25}
\begin{array}{l}
y(t;t_o,y_o,\alpha) :=\\[.25em] \left(S(t; t_0, (S_0,x_0), \alpha) - S_r(t),
\ln(x(t; t_0, (S_0,x_0), \alpha)) - \ln(x_r(t))\right)\end{array}\]
satisfies the ISS estimate (\ref{ISSest}) for all $\alpha\in
\mathcal{U}$, $t_0 \geq 0$, $t \geq t_0$, $S_0 > 0$, and $x_0 > 0$,
where $y_o=(S_0,x_0)$.
\end{theorem}

\section{Discussion on Theorem \ref{iss-track}}
\label{discuss}

Before proving the theorem, we discuss the motivations for its
assumptions, and we interpret its conclusions from both the control
theoretic and  biological viewpoints.

\begin{remark} \label{invariance}  Condition (\ref{c1ru1}) says
$(0,\infty)^2$ is positively invariant for (\ref{mod1}).   One may
also prove that  $[0,\infty)^2$ is positively invariant for
(\ref{mod1}),  as follows. Suppose the contrary. Fix $t_o\ge 0$,
$x_o\ge 0$, $S_o\ge 0$, and $\alpha\in \mathcal{U}$ for which the
corresponding trajectory $(S(t),x(t))$ for (\ref{mod1}) satisfying
$(S(t_o), x(t_o))=(S_o, x_o)$ exits $[0,\infty)^2$ in finite time.
This provides a finite constant $t_1:=\max\{\tilde t\ge t_o:
(S(t),x(t))\in [0,\infty)^2 \; \forall t\in [t_o,\tilde t]\}$.  Then
$S(t_1)=0$, since otherwise $S(t_1)>0$ and $x(t)=0$ for all $t\ge
t_1$ and then we could use the continuity of $S$ to contradict the
maximality of $t_1$. Since $\bar u<{\rm min}\{1, D_o/2\}$, the
continuity of $S$ and $x$ and the fact that $S(t_1)=0$ provide a
constant $\varepsilon>0$ such that $1+u_2(t)-S(t)\ge (1-\bar u)/2$
and $\mu(S(t))x(t)\le D_o(1-\bar u)/8$ for (almost) all $t\in [t_1,
t_1+\varepsilon]$, hence also $\dot S(t)\ge D_o(1-\bar u)/8>0$ for
all $t\in [t_1, t_1+\varepsilon]$ (since $D(t)+u_1(t)\ge D_o/2$
everywhere). Hence, $S(t)>S(t_1)=0$ for  all $t\in
[t_1,t_1+\varepsilon]$. Since $x(t)$ clearly stays in $[0,\infty)$,
this contradicts the maximality of $t_1$.  The positive invariance
of $[0,\infty)^2$ follows.
\end{remark}

\begin{remark}
\label{means} Theorem \ref{iss-track} says that in terms of the
error signals $y$, any componentwise positive trajectory of the
unperturbed chemostat dynamics (\ref{mod1}) converges to the nominal
trajectory (\ref{reftraj}), uniformly with respect to initial
conditions.  This corresponds to putting $\alpha\equiv 0$ in
(\ref{ISSest}). It also provides the additional desirable robustness
property that for an arbitrary $\mathbf{U}$-valued control function
$\alpha\in \mathcal{U}$, the trajectories of the {\em perturbed}
chemostat dynamics (\ref{mod1}) are ``not far'' from (\ref{reftraj})
for large values of time.  In other words,  they ``almost'' track
(\ref{reftraj}) with a small overflow $\gamma(|\alpha|_\infty)$ from
the ISS inequality (\ref{ISSest}). Similar results can be shown for
general choices of $x_r$ and $D(t)$. For example, we can choose any
$x_r(t)$ that admits a constant $\ell>0$ such that
\[\max\{\ell,|\dot x_r(t)|\}\; \le\; x_r(t)\; \le\; \frac{3}{4}\]
for all $t\ge 0$ and $S_r=1-x_r$. In this case, we take the dilution
rate
\[D(t)=-\frac{\dot x_r(t)}{x_r(t)}+\mu(1-x_r(t)),\] which is again
uniformly bounded above and below by positive constants. The proof
of this more general result is similar to the proof of Theorem
\ref{iss-track} we give below except with different choices of the
constants $c$ and $\kappa$.
\end{remark}

\begin{remark}
\label{mfre1} The robustness result
\begin{equation}
\label{cru1}\renewcommand{\arraystretch}{1.25}
\begin{array}{l}
|\left(S(t; t_0, (S_0,x_0), \alpha) - S_r(t), \ln(x(t; t_0,
(S_0,x_0), \alpha)) - \ln(x_r(t))\right)| \\[.25em]\le
\beta(|(S_0,x_0)|,t - t_0) + \gamma(|\alpha|_{[t_0,t]})\end{array}
\end{equation}
 of Theorem \ref{iss-track}
differs from the classical ISS condition in the following ways:

\begin{enumerate}
\item For biological reasons, negative values of the nutrient level $S$
and the species level $x$ do not make physical sense.  Hence, only
componentwise positive solutions are of interest. Therefore,
(\ref{cru1}) is not valid for all $(S_0,x_0) \in \mathbb{R}^2$ but
rather  only for $(S_0,x_0) \in (0,\infty) \times (0,\infty)$.

\item Our condition (\ref{cru1}) provides an estimate on the {\em
transformed} error  component $\ln(x(t; t_0, (S_0,x_0), \alpha)) -
\ln(x_r(t))$ instead of the more standard error $x(t; t_0,
(S_0,x_0), \alpha) - x_r(t)$.  Our reasons for using the transformed
form of the error are as follows. The function $\ln(x)$ goes to $-
\infty$ when $x$ goes to zero. This property is relevant from a
biological point of view. Indeed, in the study of biological
systems, it is important to know if the concentration of the species
is above a strictly positive constant when the time is sufficiently
large or if the concentration admits zero in its omega limit set. In
the first case, the species is called {\em persistent}. The
persistency property is frequently desirable, and it is essential to
know whether it is satisfied. Hence, the function $\ln(x(t; t_0,
(S_0,x_0), \alpha)) - \ln(x_r(t))$ has the desirable properties that
(a) it goes to $+ \infty$ if $x(t; t_0, (S_0,x_0), \alpha)$ does,
(b) it is equal to zero when $x(t; t_0, (S_0,x_0), u)$ is equal at
time $t$ to the value of $x_r$, and (c) it goes to $- \infty$ if
$x(t; t_0, (S_0,x_0), u)$ goes to zero. Therefore, roughly speaking,
if the species faces extinction, then it warns us.\end{enumerate}
\end{remark}

\begin{remark}
\label{pba} Our proof of  Theorem \ref{iss-track} is based on a
Lyapunov type analysis. Recall that a $C^1$ function
$V:\mathbb{R}^n\times [0,\infty)\to [0,\infty)$ is called an {\em
ISS Lyapunov function (ISS-LF)} for (\ref{gen}) provided there exist
$\gamma_1, \gamma_2,\gamma_3, \gamma_4\in \mathcal{K}_\infty$ such
that \begin{enumerate}
\item
$\gamma_1(|y|)\le V(y,t)\le \gamma_2(|y|)$ and \smallskip \item
$V_t(y,t)+V_y(y,t)[F(y,t)+G(y,t)u]\le
-\gamma_3(|y|)+\gamma_4(|u|)$\end{enumerate} hold for all $y\in
\mathcal{O}$, $t\ge 0$, and $u\in \mathbf{U}$. The function $V$ we
will construct in the proof of Theorem \ref{iss-track} is not an
ISS-LF for the chemostat error dynamics  because of the
specificities of the state space, which preclude the existence of
the necessary functions $\gamma_1,\gamma_2\in \mathcal{K}_\infty$ in
Condition 1 above. Hence, we cannot directly apply the  result that
the existence of an ISS Lyapunov function implies that the system is
ISS e.g. \cite[Theorem 1]{ELW00} to prove our theorem. Instead, we
prove our Theorem \ref{iss-track} directly from the decay inequality
satisfied by the time  derivative of $V$ along the trajectories. The
proof that the decay inequality implies ISS is very similar to that
part of the proof of \cite[Theorem 1]{ELW00}, so we only sketch that
part of our proof in the appendix.
\end{remark}

\begin{remark}
\label{bound-u} Our estimate (\ref{cru1}) would not hold if we had
instead chosen the full control set $\mathbf{U}=\mathbb{R}^2$. In
fact, taking the disturbance $\alpha\equiv (u_1, u_2)=(0,-1)$ and
any initial condition $(S(t_o),x(t_o))=(S_0,x_o)\in (0,\infty)^2$ in
(\ref{mod1}) would give $S(t; t_0, (S_0,x_0), \alpha) \to 0$ and so
also $\ln(x(t; t_0, (S_0,x_0), \alpha))) \to - \infty$ as $t\to
+\infty$ (since $|u_1(t)|\le \bar u<D_o/2\le D(t)/2$ almost
everywhere).
 Therefore extinction would occur and (\ref{cru1}) would not
be satisfied.  On the other hand, if our  set $\mathbf{U}$ is
replaced by $\mathbf{U}^\sharp:=[-\bar u,+\bar u]^2$ for any fixed
constant $\bar u\in (0,\min\{1,D_o\})$, then the chemostat error
dynamics instead satisfies the less stringent {\em integral} ISS
property; see Remark \ref{iISS-rk} below for details.
\end{remark}

\section{Proof of Theorem \ref{iss-track}}\label{main}

The proof of (\ref{c1ru1}) is immediate from the structure of the
dynamics (\ref{mod1}) and the fact that $\bar u<1$ (which imply that
 $\dot S\ge 0$ when $S>0$ is
sufficiently small);  see Remark \ref{invariance} for a similar
argument.  It remains to prove the ISS estimate (\ref{cru1}) for
suitable functions $\beta\in \mathcal{KL}$ and $\gamma\in
\mathcal{K}_\infty$.

Throughout the proof, all (in)equalities should be understood to
hold globally unless otherwise indicated.  Also,  we repeatedly use
the simple ``(generalized) triangle inequality'' relation
\begin{equation}
\label{triangle} pq\; \; \le\; \;  d p^2+\frac{1}{4d}q^2
\end{equation}
for various choices of $p\ge 0$, $q\ge 0$, and $d>0$ that we specify
later.

Fix $t_o\ge 0$, $S_o>0$, $x_o>0$, and $\alpha\in \mathcal{U}$, and
let $[t_o,\infty)\ni t\mapsto (S(t),x(t))$ denote the corresponding
solution of (\ref{mod1}) satisfying $(S(t_0),x(t_o))=(S_o,x_o)$. For
simplicity, we write $\alpha(t)$ as $(u_1,u_2)$, omitting the time
argument as before. We first write
 the error equation for the  variables
\begin{equation}
\label{impo} (\tilde z, \tilde \xi) = (z - z_r, \xi - \xi_r)
\end{equation}
where $\xi = \ln x$, $z = S + x$, $z_r(t) = S_r(t) + x_r(t)=1$, and
$\xi_r(t) = \ln x_r(t)$. One easily checks that \[
\begin{array}{rcl}
\dot z(t)&=& [D(t)+u_1(t)][1+u_2(t)-z(t)]\\
\dot x(t)&=& x(t)[\mu(z(t)-x(t))-D(t)-u_1(t)]\, .
\end{array}
\]
Therefore,  since $z_r\equiv 1$ (which implies $\dot z_r(t)=
[D(t)+u_1(t)][1-z_r(t)]$), our formula (\ref{mbound}) for $\mu$
immediately implies that the (transformed) error  \[(\tilde z,
\tilde \xi)(t)=(z(t)-z_r(t), \ln x(t)-\ln x_r(t))\] satisfies the
{\em chemostat error dynamics}
\begin{equation}
\label{dis} \left\{
\begin{array}{rcl}
\dot{\tilde{z}} & = & - [D(t) + u_1(t)] [\tilde z-u_2(t)]
\\[.5em]
\dot{\tilde{\xi}} & = & ma\dfrac{\tilde{z} -
e^{\xi_r(t)}(e^{\tilde{\xi}} - 1)}{(a + z - e^{\xi})(a + z_r(t) -
e^{\xi_r(t)})} - u_1(t).
\end{array}
\right.
\end{equation}
We are going to show that (\ref{dis}) has the Lyapunov function
\begin{equation}
\label{VChoice}
\begin{array}{l}
 V(\tilde \xi, \tilde z) =
e^{L_3(\tilde{\xi}, \tilde z)} - 1,\; \; \; {\rm where}\; \;
L_3(\tilde \xi, \tilde z)=L_1(\tilde \xi)+\kappa L_2(\tilde z),\\
L_1(\tilde \xi)=e^{\tilde \xi}-\tilde \xi-1,\; \; {\rm and}\; \;
L_2(\tilde z) = \frac{1}{D_o - \bar u} \tilde z^2
\end{array}
\end{equation}
and $\kappa>0$ is the constant defined in (\ref{ck}). {}From the
explicit expressions $z_r(t) = 1$ and  $e^{\xi_r(t)} = \frac{1}{2} +
\frac{1}{4}\cos(t)\ge \frac{1}{4}$, we deduce that the time
derivative of $L_1$ along trajectories of (\ref{dis})  satisfies
\[\renewcommand{\arraystretch}{2.5}
\begin{array}{rcl}
\dot{L}_1 & = & ma\dfrac{(e^{\tilde{\xi}} - 1)\tilde{z} -
e^{\xi_r(t)} (e^{\tilde{\xi}} - 1)^2}{(a + z - e^{\xi})(a + z_r(t) -
e^{\xi_r(t)})} - (e^{\tilde{\xi}} - 1) u_1(t)\\
& \le & ma\dfrac{ - \frac{1}{4 a + 2 - \cos(t)}(e^{\tilde{\xi}} -
1)^2 + \frac{4}{4 a + 2 - \cos(t)}|e^{\tilde{\xi}} -
1||\tilde{z}|}{(a + z - e^{\xi})} - (e^{\tilde{\xi}} - 1) u_1(t)
\\
& \leq & \dfrac{ - \frac{ma}{4 a + 3}(e^{\tilde{\xi}} - 1)^2 +
\frac{4 ma}{4 a + 1}|e^{\tilde{\xi}} - 1||\tilde{z}|}{(a + z -
e^{\xi})} - (e^{\tilde{\xi}} - 1) u_1(t) ,
\end{array}
\]
where we also used the fact that $z-e^\xi=S\ge 0$.  Since
$D(t)+u_1(t)\ge D_o-\bar u$ everywhere, one readily checks that
along the trajectories of (\ref{dis}),
\begin{equation}
\label{s1mple}
\begin{array}{rclcl}
\dot{L}_2  &\leq & - \tilde z^2 + \tilde c|\tilde z| |u_2(t)|&\le
&-\frac{1}{2}\tilde z^2+cu^2_2(t)
\end{array}
\end{equation}
where $\tilde c = 2(\bar D + \bar u)(D_o - \bar u)^{-1}$, the
constant $c$ is defined by (\ref{ck}), and the last inequality used
(\ref{triangle}) with the choices $p=|\tilde z|$, $q=|u_2(t)|$, and
$d=1/(2\tilde c)$.  The fact that $\frac{1}{2}\tilde c^2\le c$
easily follows because $\bar u\le \frac{1}{2}D_o$. Along the
trajectories of (\ref{dis}),
\begin{equation}
\label{s2mple}
\begin{array}{rcl}
\dot{L}_3 & \leq & \dfrac{ - \frac{ma}{4 a + 3}(e^{\tilde{\xi}} -
1)^2 + \frac{4 ma}{4 a + 1}|e^{\tilde{\xi}} - 1||\tilde{z}|}{(a + z
- e^{\xi})} - (e^{\tilde{\xi}} - 1) u_1(t)\\&& - \dfrac{1}{2}\kappa
\tilde z^2 + \kappa c u^2_2(t).
\end{array}
\end{equation}
We distinguish between two cases.

\noindent {\em \underline{Case 1a}: $z(t) \leq 2$}. Then since
$z-e^\xi=S \ge 0$, we get
\begin{equation}
\label{s2h}\renewcommand{\arraystretch}{1.5}
\begin{array}{rcl}
\dot{L}_3 & \leq & - \dfrac{m a}{(4a + 3)(a + 2)}(e^{\tilde{\xi}} -
1)^2 + \dfrac{4 m}{4a + 1} |e^{\tilde{\xi}} - 1||\tilde{z}|\\&& -
(e^{\tilde{\xi}} - 1) u_1(t) - \frac{1}{2}\kappa \tilde z^2 + \kappa
c u^2_2(t).
\end{array}
\end{equation}
Using the triangle inequality (\ref{triangle}) with the choices
\[
p=|e^{\tilde \xi}-1|,\; \; \;  q=|\tilde z|, \; \; \; {\rm and}\; \;
\;  d=\frac{a(4a+1)}{8(4a+3)(a+2)},\] we deduce from (\ref{s2h})
that
\begin{equation}
\label{s21}\renewcommand{\arraystretch}{1.5}
\begin{array}{rcl}
\! \! \dot{L}_3 & \leq & - \dfrac{m a}{2(4a + 3)(a +
2)}(e^{\tilde{\xi}} - 1)^2\\[1em]&& - \left[\dfrac{\kappa}{2} -
\dfrac{8(4a + 3)(a + 2)m}{a (4 a + 1)^2}\right] \tilde z^2 -
(e^{\tilde{\xi}} - 1) u_1 + \kappa c u^2_2.
\end{array}
\end{equation}
 \noindent {\em \underline{Case 2a}: $z(t) \geq 2$}. Since $z(t) - e^{\xi(t)}
\geq 0$, it follows that $e^{- \xi_r(t)} z(t) - e^{\tilde \xi(t)}
\geq 0$. Therefore, since $x_r\ge  1/4$,
\begin{equation}
\label{a4h}\begin{array}{l}
 e^{\tilde \xi(t)}  \; \leq\;   e^{- \xi_r(t)} z(t)  \;
\leq\;   e^{-\ln(1/4)}z(t)= 4 z(t),\\ {\rm hence}\; \; \; - 1 \;
\leq \;  e^{\tilde \xi(t)} - 1  \; \leq  \; 4 z(t) - 1.
\end{array}\end{equation} Since $z(t) \geq 2$ and $z_r=1$, we have  $\tilde
z(t) \geq 1$. As $\tilde z=z-z_r=z-1$, condition (\ref{a4h}) gives
\begin{equation} - \tilde z(t) \;  \leq\;  e^{\tilde \xi(t)} - 1
\; \leq\; 3 + 4\tilde z(t) \; \leq\;  7 \tilde z(t), \; \; \; {\rm
hence}\; \; \; |e^{\tilde \xi(t)} - 1| \leq 7 \tilde z(t).
\end{equation}
{}From this last inequality and the inequality $z - e^{\xi} \geq 0$,
we deduce from dropping the first term in (\ref{s2mple}) that
\begin{equation}
\label{u2upue}
\begin{array}{rcl}
\dot{L}_3 & \leq & \dfrac{28 m}{4 a + 1} \tilde z^2 -
(e^{\tilde{\xi}} - 1) u_1 - \dfrac{\kappa}{2} \tilde z^2 + \kappa c
u^2_2
\\[1em]
& \leq & - \dfrac{1}{200}\kappa (e^{\tilde \xi} - 1)^2 -
\left[\dfrac{1}{4} \kappa - \dfrac{28 m}{4 a + 1}\right] \tilde z^2
- (e^{\tilde{\xi}} - 1) u_1 + \kappa c u^2_2 .
\end{array}
\end{equation}

\noindent We deduce from our choice (\ref{ck}) of $\kappa$,
(\ref{s21}), and (\ref{u2upue}) that in Cases 1a-2a,
\begin{equation}
\label{a2upue}
\begin{array}{rcl}
\dot{L}_3 & \leq & - C_1 [(e^{\tilde \xi} - 1)^2 + \tilde z^2] -
(e^{\tilde{\xi}} - 1) u_1 + \kappa c u^2_2
\end{array}
\end{equation}
where $C_1$ is defined in (\ref{ck}).

Using (\ref{triangle}) with $p=|{\rm exp}(\tilde \xi(t))-1|$,
$q=|u_1(t)|$, and $d=\frac{C_1}{2}$, and then the upper bounds of
$u_1$ and $u_2$, we deduce from (\ref{a2upue}) that
\begin{equation}
\label{a2uiue}\renewcommand{\arraystretch}{2.5}
\begin{array}{rcl}
\! \! \! \!  \dot{L}_3 & \leq & - \dfrac{C_1}{2} [(e^{\tilde \xi} -
1)^2 + \tilde z^2] +
 \dfrac{1}{2 C_1} \bar u |u_1| + \kappa c \bar u |u_2|
\\
 & \leq & - \dfrac{C_1}{2} [(e^{\tilde \xi} - 1)^2 + \tilde z^2] + C_2
 |u|,\; \; \;  {\rm where} \; \; \; C_2 := \left(\dfrac{1}{C_1}  + 2\kappa c\right) \bar
u
\end{array}
\end{equation}
and where the last inequality used the relationship $|u|_1\le
2|u|_2$ between the $1$-norm and the $2$-norm. We consider  two
additional cases.

\noindent \underline{\em Case 1b}: $(e^{\tilde \xi(t)} - 1)^2 +
\tilde z^2(t) \geq \frac{1}{2}$. Then (\ref{a2uiue}) gives
\begin{equation}
\label{s2dide}
\begin{array}{rcl}
\dot{V} & \leq & e^{L_3(\tilde{\xi},\tilde{z})}\left(-
\dfrac{C_1}{4} + C_2 |u(t)|\right).
\end{array}
\end{equation}
Next notice that
 (\ref{constraint}) and our choice of $C_1\in (0,1]$ give
 \begin{equation}
\label{s2dodo} \bar u \; \leq\; \frac{C_1}{8 C_2},\; \; \; \; {\rm
hence} \; \; \; \; \dot{V}
 \; \leq\;  - \dfrac{C_1}{8}e^{L_3(\tilde{\xi},\tilde{z})} \; \leq\;  -
\dfrac{C_1}{8} V(\tilde \xi, \tilde z).
\end{equation}

\noindent \underline{\em Case 2b}: $(e^{\tilde \xi(t)} - 1)^2 +
\tilde z^2(t) \leq \frac{1}{2}$. Then $(\tilde \xi(t),\tilde z(t))$
is in a suitable bounded set, so since
\[\tilde F(L)=(e^L-L-1)(e^L-1)^{-2},\;\; \;
{\tilde G}(L)= (e^L-1)L^{-1}
\] are locally bounded when defined to be zero at $L=0$, one can
readily use (\ref{a2uiue}) to compute constants $C_3 > 0$ and $C_4
> 0$ such that
\begin{equation}\label{een}
\begin{array}{rclcl} \dot{V} & \leq & - C_3 L_3(\tilde{\xi},
\tilde z) + C_2 |u(t)| &\le & -C_4 V(\tilde{\xi}, \tilde z)+ C_2
|u(t)|
\end{array}
\end{equation}
where $\tilde F$ was used to get $C_3$ and $\tilde G$ was used to
get $C_4$. It follows from (\ref{s2dodo})-(\ref{een}) that, in Cases
1b-2b,
\begin{equation}
\label{goal}
\begin{array}{rcl}
\dot{V} & \leq & - C_5 V(\tilde{\xi}, \tilde z) + C_2 |u(t)|
\end{array}
\end{equation}
with $C_5 = \inf\left\{C_4, \frac{C_1}{8}\right\}$.  Condition
(\ref{goal}) is the classical ISS Lyapunov function decay condition
for the transformed error dynamics evolving on our restricted state
space. Therefore, a slight variant of the classical ISS arguments
combined with (\ref{goal}) give the ISS estimate asserted by Theorem
\ref{iss-track}.  For details, see the appendix below. This
concludes the proof.

\begin{remark}
\label{iISS-rk} If our control set $\mathbf{U}$ is replaced by the
larger control set $\mathbf{U}^\sharp:=[-\bar u,+\bar u]^2$ for any
fixed constant $\bar u\in (0,\min\{1,D_o\})$, then the error
dynamics (\ref{dis}) instead satisfies the less stringent {\em
integral} ISS property.  The relevant definitions are as follows. We
say that (\ref{gen}) is {\em integral input-to-state stable (iISS)}
provided there exist $\delta_1,\delta_2\in \mathcal{K}_\infty$ and
$\beta\in \mathcal{KL}$ such that \[ \tag{iISS} \delta_1(|y(t; t_o,
y_o, \alpha)|)\; \le\; \beta(|y_o|,
t-t_o)+\int_{t_o}^{t+t_o}\delta_2(|\alpha(r)|)dr\] everywhere for
all measurable essentially bounded functions $\alpha:[0,\infty)\to
\mathbf{U}^\sharp$.  This condition is less restrictive than ISS
since e.g. $\dot y=-\arctan(y)+u$ is iISS but not ISS \cite{ASW00}.
An {\em iISS-LF} for (\ref{gen}) (with controls in
$\mathbf{U}^\sharp$) is then defined to be a $C^1$ function
$V:\mathbb{R}^n\times [0,\infty)\to [0,\infty)$ for which there are
$\gamma_1, \gamma_2, \gamma_4\in \mathcal{K}_\infty$
 and a positive definite function $\gamma_3$ (i.e.,
 $\gamma_3:[0,\infty)\to[0,\infty)$ is continuous and zero only at
 zero)
such that Conditions 1-2 in Remark \ref{pba} hold everywhere. This
is less restrictive than the ISS-LF condition since $\gamma_3$ need
not be of class $\mathcal{K}_\infty$. Arguing as in the proof of
Theorem \ref{iss-track} up through (\ref{a2uiue}) and solving the
appropriate constrained minimum problem to get $\gamma_3$ shows that
$V=L_3$ satisfies the iISS Lyapunov function decay condition
(namely, Condition 2 from Remark \ref{pba})
 for the error dynamics (\ref{dis}) and the control set
$\mathbf{U}^\sharp$ using $\gamma_3(s)=C_1(e^{-s}-1)^2/2$.
Therefore, this system is in fact iISS, by a slight variant of the
proof of the iISS estimate in \cite[Theorem 1]{ASW00}.  We leave the
details to the reader.
\end{remark}

\section{Stability in the Presence of Several Species}
\label{several} Theorem \ref{iss-track} shows that the stability of
the reference  trajectory (\ref{reftraj}) is robust with respect to
small perturbations of the dilution rate and initial concentration.
To further demonstrate the robustness of our results, we
 next show
that the stability of (\ref{reftraj}) is also maintained when the
model (\ref{model}) is augmented to include additional species that
are being driven to extinction, in the following sense.

We assume for simplicity that $u_1\equiv u_2\equiv 0$.
 Consider the augmented system
 \begin{equation}\renewcommand{\arraystretch}{.8}
\label{nmodel} \! \! \! \! \! \left\{\begin{array}{rcl} {\dot S} & =
& D(t)(1 - S) - \mu(S) x - \displaystyle \sum_{i=1}^n\nu_i(S)y_i
\\
\dot{x} & = & x (\mu(S) - D(t))
\\[.5em]
{\dot y_i} & = & y_i(\nu_i(S) - D(t)),\;\; i=1,\dots, n ,
\end{array}\right.
\end{equation}

\noindent where $\mu$ is as in (\ref{mbound}) and $\nu_i$ is
continuous and increasing and satisfies $\nu_i(0) = 0$ for
$i=1,2,\ldots, n$. The variables $y_i$ represent the levels of $n$
additional species. We choose $D$ and $D_o$ as in (\ref{chli}) and
(\ref{aer1}), and we assume  $\nu_i(1)<D_o$ for $i=1,2,\ldots, n$.
(This assumption is, in a sense, natural because one can easily
check that it ensures that each species concentration $y_i$
converges to zero. Indeed, the fact that, for all $t \geq 0$, $D(t)
> D_o > 0$ and $\mu(S(t)) x(t) +
\sum_{i=1}^n\nu_i(S(t))y_i(t) \geq 0$ ensures, in combination with
the inequalities $\nu_i(1)<D_o$, that there exists an instant $T >
0$ and a constant $c > 0$ such that for all $t \geq T$ and for
$i=1,2,\ldots, n$, $\dot{y}_i(t) \leq - c y_i(t)$. This implies that
$y_i$'s converge to $0$ exponentially.) We show that the transformed
error
\begin{equation}
\label{nte} (\tilde z,\tilde \xi,\tilde y)\; :=\;
(S+x-S_r-x_r,\ln(x)-\ln(x_r),y)
\end{equation}
between any componentwise positive solution $(S,x,y)$ of
(\ref{nmodel}) and the  reference trajectory
 $(S_r,x_r,0, \ldots, 0) =
\left(\frac{1}{2} - \frac{1}{4}\cos(t) ,\frac{1}{2} +
\frac{1}{4}\cos(t), 0,\ldots,0\right)$ converges exponentially to
the zero vector as $t\to +\infty$.

To this end, notice that in the coordinates (\ref{nte}), the system
(\ref{nmodel}) becomes
\begin{equation}
\label{dmodel}\renewcommand{\arraystretch}{2.25} \! \! \! \! \!
\left\{\begin{array}{rcl} \dot{\tilde z} & = & - D(t) \tilde z(t) -
\displaystyle \sum_{i=1}^n\nu_i(S)y_i
\\
\dot{\tilde \xi} & = &  ma\dfrac{\tilde{z} - e^{\xi_r(t)}
(e^{\tilde{\xi}} - 1)}{(a + z - e^{\xi})(a + z_r(t) - e^{\xi_r(t)})}
\\
\dot{y}_i & = & y_i(\nu_i(S) - D(t)),\;\; i=1,\dots, n ,
\end{array}\right.
\end{equation}
by the same calculations that led to (\ref{dis}). Set $L_1(\tilde
\xi)={\rm exp}(\tilde \xi)-\tilde \xi-1$ where ${\rm exp}(r):=e^r$.
 Since $e^{\xi_r}\ge 1/4$,
$z_r\equiv 1$, $0\le z_r-e^{\xi_r}=S_r=1/2-\cos(t)/4\leq 1$, and
$0\le z-e^\xi=\tilde z+1-e^\xi\le 2+\tilde z^2$,
 we deduce that the
derivative of \[L_3(\tilde \xi,\tilde z):=L_1(\tilde
\xi)+\frac{4m}{aD_o}\tilde z^2\] along the trajectories of
(\ref{dmodel}) satisfies
\begin{equation}
\label{dmo1}\renewcommand{\arraystretch}{2.25}
\begin{array}{rcl}
\dot{L}_3 & \leq & ma\dfrac{\tilde z(e^{\tilde
\xi}-1)-\frac{1}{4}(e^{\tilde
\xi}-1)^2}{(a+z-e^\xi)(a+z_r-e^{\xi_r})}-\dfrac{8m}{a}\tilde z^2-
\dfrac{8m}{aD_o}\tilde{z} \displaystyle \sum_{i=1}^n\nu_i(S)y_i\\
&\le & ma\dfrac{4\tilde z^2-\frac{1}{16}(e^{\tilde
\xi}-1)^2}{(a+z-e^\xi)(a+z_r-e^{\xi_r})}-\dfrac{8m}{a}\tilde z^2-
\dfrac{8m}{aD_o}\tilde{z} \displaystyle \sum_{i=1}^n\nu_i(S)y_i\\
&\leq& - \dfrac{ma(e^{\tilde{\xi}} - 1)^2}{16(a+1) (a + 2 +
\tilde{z}^2)} - \dfrac{4m}{a}\tilde{z}^2\\[.5em]
&&-2\left[\dfrac{\sqrt{m}}{\sqrt{a}}{\tilde
z}\right]\left[\dfrac{4\sqrt{m}}{\sqrt{a}{D_o}} \displaystyle
\sum_{i=1}^n\nu_i(S)y_i\right]
\\[1em]
& \leq & - \dfrac{ma(e^{\tilde{\xi}} - 1)^2}{16(a+1) (a + 2 +
\tilde{z}^2)} - \dfrac{3m}{a}\tilde{z}^2 + \dfrac{16 m}{ a
D^2_o}\left[\displaystyle \sum_{i=1}^n\nu_i(S)y_i\right]^2
\end{array}
\end{equation}
where the second inequality is by (\ref{triangle}) with $p=|\tilde
z|$, $q=|e^{\tilde \xi}-1|$, and $d=4$ and the last inequality used
 the relation $J^2+K^2\ge -2JK$ for real values $J$
and $K$. On the other hand, since $\nu_i(1)< D_o$ for each $i$, the
form of the dynamics for $S$ and the nonnegativity of $\mu$ and the
$\nu_i$'s along our componentwise positive trajectories imply that
there exist $\varepsilon > 0$ and $T \geq 0$ such that (i)  $S(t)
\leq 1 + \varepsilon$  for  all  $t\ge T$ and (ii)  $\nu_i(1 +
\varepsilon) < D_o$  for all $i=1,2,\ldots, n$. We deduce that, for
all $i=1,2,\ldots, n$ and for all $t \geq T$,
\begin{equation}
\label{exs1}\renewcommand{\arraystretch}{2}
\begin{array}{rclclcl}
\! \! \! \dfrac{1}{2}\dfrac{d}{dt}y^2_i &\leq& (\nu_i(S(t)) - D(t))
y_i^2 &\leq& (\nu_i(1 + \varepsilon) - D_o) y_i^2 &\leq& - \delta
y_i^2,
\end{array}\end{equation}
where $\delta=D_o-\max\{\nu_i(1 + \varepsilon): i=1,2,\ldots,
n\}>0$. Hence,  each $y_i(t)$ converges exponentially to zero.

Next notice that along each pair $(\tilde \xi(t),\tilde z(t))$, the
function \[\Delta(\tilde \xi,\tilde z):= \dfrac{ma(e^{\tilde{\xi}} -
1)^2}{16(a+1) (a + 2 + \tilde{z}^2)}
\]
is positive if and only if $\tilde\xi \ne 0$. By (\ref{dmo1}) and
(\ref{exs1}),  the time derivative of
\begin{equation}
\label{dmo3} L_4(\tilde{z}, \tilde{\xi}, y_1, ..., y_n) =
L_3(\tilde{z}, \tilde{\xi}) + A \displaystyle\displaystyle \sum_{i =
1}^{n} y_i^2, \; \; \; {\rm where}\; \; A:=\frac{16mn^2}{a\delta}
\end{equation}
along the trajectories of (\ref{dmodel}) satisfies
\begin{equation}\label{surmise}\renewcommand{\arraystretch}{2.25}
\begin{array}{rcl}
\dot{L}_4 &\le& - \Delta(\tilde \xi,\tilde z) -
\dfrac{3m}{a}\tilde{z}^2
 +\dfrac{16m}{aD^2_o}\left[\displaystyle
\sum_{i=1}^n\nu_i(S)y_i\right]^2-2A\delta\displaystyle
\sum_{i=1}^ny^2_i\\
&\le& - \Delta(\tilde \xi,\tilde z) - \dfrac{3m}{a}\tilde{z}^2
 +\dfrac{16mn^2}{aD^2_o}\displaystyle
\sum_{i=1}^n\nu^2_i(S)y^2_i-2A\delta\displaystyle
\sum_{i=1}^ny^2_i\\&\le&  - \Delta(\tilde \xi,\tilde z)
-\dfrac{3m}{a}\tilde z^2-\dfrac{16mn^2}{a}\displaystyle
\sum_{i=1}^ny^2_i\; \;  =:\; \;  -M(\tilde \xi,\tilde z,y)
\end{array}\end{equation}

\noindent provided $t>T$ where $T$ is chosen as above.  (The second
inequality in (\ref{surmise}) follows because for any nonnegative
 $a_k$, we get $a_k\le (\sum_{i=1}^na^2_i)^{1/2}$ which we sum
and then square to get $(\sum_{i=1}^na_i)^2\le n^2
\sum_{i=1}^na^2_i$.  The last inequality used  $\nu_i(S(t))\le
\nu_i(1+\varepsilon)<D_o$ for all $t\ge T$ and our choice of $A$. )

It is tempting to surmise from (\ref{surmise}) and the structure of
$L_4$ that $L_4$ is a Lyapunov function for (\ref{dmodel}) since
then we could use standard Lyapunov function theory to conclude that
$(\tilde \xi(t),\tilde z(t), y(t))$ asymptotically converges to
zero. However, such an argument would not be technically correct,
since the state space of (\ref{dmodel}) is not ${\mathbb R}^{n+2}$
(because the original augmented chemostat model (\ref{nmodel})
  is
only defined for componentwise nonnegative values of the state).
Instead, we argue as follows (in which we may assume  for simplicity
that the initial time for the augmented error dynamics is zero).

For any $t \geq 0$, integrating the last inequality of
(\ref{surmise}) over $[0,t]$ gives
\begin{equation}
\label{da2}
\begin{array}{l}
 L_4\left(\tilde{z}(t), \tilde{\xi}(t), y(t)\right) - L_4\left(\tilde{z}(0),
 \tilde{\xi}(0), y(0)\right)\\\leq
 - \displaystyle\int_{0}^{t} M\left(\tilde \xi(l),\tilde
z(l),y(l)\right) dl\; .
\end{array}
\end{equation}
It follows that, for all $t \geq 0$,
\begin{equation}
\label{da3}
\begin{array}{rcl}
L_1\left(\tilde \xi(t)\right) &\leq&  L_4\left(\tilde{z}(0),
\tilde{\xi}(0), y(0)\right)\; .
\end{array}
\end{equation}
Therefore $\xi(t)=\tilde \xi(t)+\xi_r(t)$ is a bounded function.
Similarly, $\tilde{z}(t)$ and $y(t)$ are bounded. We deduce that
$\tilde\xi$, $\tilde z$, and the components of $y$ are uniformly
continuous, since their time derivatives (\ref{dmodel}) are bounded.
Reapplying (\ref{da2}) therefore implies \[\int_{0}^{+ \infty}
M(\tilde \xi(l),\tilde z(l),y(l))  dl\] is finite. It follows from
Barbalat's lemma \cite[p.323]{K02} and the structure of the function
$M$ that $(\tilde \xi(t), \tilde z(t), y(t))\to 0$ as $t\to+\infty$.
This establishes our stability condition for the multi-species
model.

\begin{remark}
 Notice that $-\dot L_4$ is bounded below by a quadratic of the
form $\bar c|(\tilde \xi, \tilde z,y)|^2$ along the trajectories of
(\ref{dmodel}), and that $L_4$ is bounded above and below by such
quadratics along the trajectories as well,  since the trajectories
are bounded. {}From this fact  and (\ref{surmise}), one can deduce
that the trajectories $(\tilde \xi(t), \tilde z(t),y(t))$ converge
{\em exponentially} to zero.
\end{remark}

\section{Simulation}
\label{simulations} 
\begin{figure}[t]
\begin{center}
\scalebox{.5}{\includegraphics{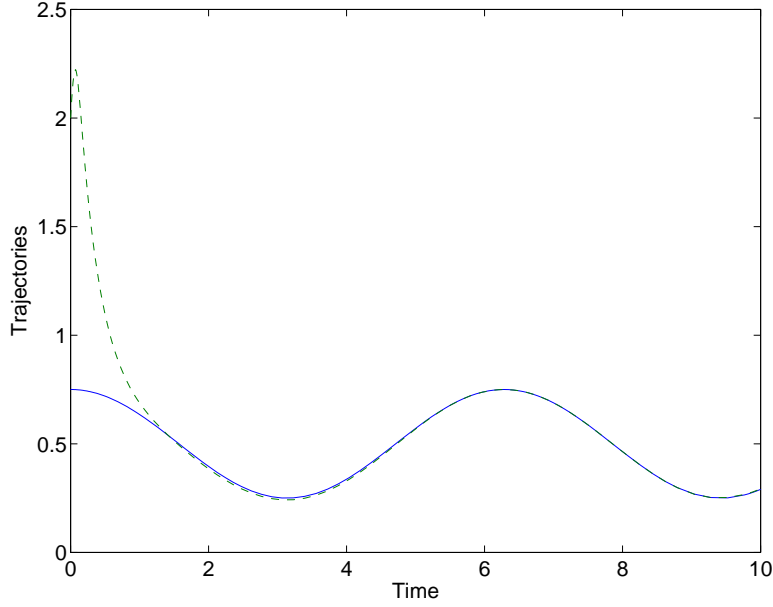}}
\end{center}
\caption{State Trajectory Component $x(t)$ (dashed) and Reference
Trajectory $x_r(t)$ (solid) for Chemostat \label{xsol}}
\end{figure}

To validate our convergence result, we simulated
the dynamics (\ref{mod1}) with the initial values $x(0)=2$ and
$S(0)=1$ and the reference trajectory $x_r(t)$, using the parameters
$m=10$ and $a=\frac{1}{2}$ and $t_o=0$.  In this case, the lower
bound on $D(t)$ provided by (\ref{aer1}) is $D_o=7/3$.  It follows
from Remark \ref{iISS-rk} that the convergence of $x(t)$ to $x_r(t)$
is robust to disturbances that are valued in $[-\bar u, \bar u]^2$
for any positive constant $\bar u<\min\{1,D_o\}=1$, in the sense of
integral input-to-state stability.  Moreover, using the estimate
(\ref{a2uiue}), one easily checks that in this case, estimate (iISS)
on p.\pageref{iISS-rk} above holds with $\delta_2(r)=2C_2r$; cf. the
proof of \cite[Theorem 1]{ASW00}.
 For our simulation, we took the
disturbance  $u_1(t)=0.5e^{-t}$ on the dilution rate and
$u_2(t)\equiv 0$. This gave the plot of $x(t)$ and $x_r(t)$ against
time in Figure \ref{xsol}. Our simulation shows that the state
trajectory $x(t)$ closely tracks the reference trajectory $x_r(t)$
even in the presence of small disturbances and so validates our
findings.

\section{Conclusions}
\label{concl}

The chemostat model is a useful framework for modeling species
competing for nutrients.  For the case of
 one species competing for one
nutrient and a suitable time-varying dilution rate, we proved
stability of an appropriate reference trajectory. Moreover, we found
that the stability was maintained even if the model is augmented
with other species that are being driven to extinction, or if there
are disturbances of appropriately small magnitude acting on the
dilution rate and input nutrient concentration.

\smallskip

\section*{Appendix}
For completeness, we provide the slight variant of the classical ISS
arguments needed to finish the proof of Theorem \ref{iss-track}.
Multiplying through (\ref{goal}) by $e^{C_5l}$ and applying the
standard ``variation of parameters'' formula to $[t_o,t]\ni l\mapsto
V(\tilde \xi(l),\tilde z(l))$ (by integrating between $t_0\ge 0$ and
$t \geq t_o$) gives
\begin{equation}
\begin{array}{rcl}
V(\tilde{\xi}(t), \tilde z(t)) &\leq& e^{(t_0 - t) C_5}
V(\tilde{\xi}(t_0), \tilde z(t_0))+ C_2 |u|_{[t_o,t]} \; ,
\end{array}
\end{equation}
where we enlarged $C_2$ without relabeling. We deduce that
\[
L_1(\tilde{\xi}(t)) + \frac{\kappa}{D_o - \bar u} \tilde z^2(t) \leq
\ln\left(1 + e^{(t_0 - t) C_5} V(\tilde{\xi}(t_0), \tilde z(t_0)) +
C_2 |u|_{[t_o,t]}\right)\] where $L_1$ is defined in
(\ref{VChoice}). Since $e^r - 1 - r \geq \frac{1}{2}r^2$ and $\ln(1
+ r) \leq r$ for all $r \geq 0$, we deduce from the formula for $V$
that
\begin{equation}\label{49}\begin{array}{rcl}
\dfrac{1}{2}\tilde{\xi}^2(t) +  \dfrac{\kappa}{D_o - \bar u} \tilde
z^2(t) & \leq&   e^{(t_0 - t) C_5}
\Omega\left(|(\tilde{\xi}(t_0),\tilde z(t_0))|\right) + C_2
|u|_{[t_o,t]}\\&& {\rm where}\; \; \; \; \Omega(r) = e^{e^{r} - 1 -
r + \frac{\kappa}{D_o - \bar u} r^2} - 1\; .
\end{array}
\end{equation}
In particular, $\Omega\in{\mathcal K}_{\infty}$.  {}From (\ref{49})
and the inequality $\sqrt{a+b}\le \sqrt{a}+\sqrt{b}$,
\begin{equation}
\label{tal}
\begin{array}{rcl}
|\tilde{\xi}(t)| & \leq & \sqrt{2 e^{(t_0 - t) C_5}
\Omega\left(|(\tilde{\xi}(t_0),\tilde z(t_0))|\right)} + \sqrt{2 C_2
|u|_{[t_o,t]}} \; \; \; \; {\rm and}
\end{array}
\end{equation}
\begin{equation}
\begin{array}{rcl}
|\tilde z(t)| & \leq & \sqrt{e^{(t_0 - t) C_5} \frac{D_o - \bar
u}{\kappa} \Omega\left(|(\tilde{\xi}(t_0),\tilde z(t_0))|\right)} +
\sqrt{C_2 \frac{D_o - \bar u}{\kappa}|u|_{[t_o,t]}}
\end{array}\; .
\end{equation}
The relations  $\tilde z=S-S_r+e^{\xi_r}(e^{\tilde \xi}-1)$ and
$e^{a+b}-1\le \frac{1}{2}(e^{2a}-1)+\frac{1}{2}(e^{2b}-1)$ give
\[\renewcommand{\arraystretch}{2}
\begin{array}{rcl}
|S(t) - S_r(t)| & \leq & |\tilde z(t)| +  (e^{|\tilde{\xi}(t)|} -
1)\\
&  \hspace{-.75in}\leq & \hspace{-.35in} \sqrt{e^{(t_0 - t) C_5}
\frac{D_o - \bar u}{\kappa} \Omega\left(|(\tilde{\xi}(t_0),\tilde
z(t_0))|\right)} + \sqrt{C_2 \frac{D_o - \bar
u}{\kappa}|u|_{[t_o,t]}}
\\
& & \hspace{-.5in}+ \frac{1}{2} \left(e^{2\sqrt{2 e^{(t_0 - t) C_5}
\Omega\left(|(\tilde{\xi}(t_0),\tilde z(t_0))|\right)}} - 1\right) +
\frac{1}{2} \left(e^{2\sqrt{2 C_2 |u|_{[t_o,t]}}} - 1\right) .
\end{array}
\]

The desired ISS estimate (\ref{ISSest}) now follows  immediately
from this last inequality and (\ref{tal}) with the choices
\[\renewcommand{\arraystretch}{1.5}\begin{array}{l}\renewcommand{\arraystretch}{1.25}
\beta(s,t)=4\sqrt{\Omega(s){\rm exp}(-C_5t)\left\{1+\frac{D_o-\bar
u}{\kappa}\right\}}+ {\rm exp}\left(4\sqrt{\Omega(s){\rm
exp}(-C_5t)}\right)-1\\
{\rm and}\; \; \gamma(r)=4\sqrt{C_2\left(1+(D_o-\bar
u)/\kappa\right)r}+{\rm exp}(4\sqrt{c_2r})-1.\end{array}
\]
This completes the proof of Theorem \ref{iss-track}.

\section*{Acknowledgments}

Part of this work was done while P. De Leenheer and F. Mazenc
visited Louisiana State University (LSU). They thank LSU for the
kind hospitality they enjoyed during this period. F. Mazenc thanks
Claude Lobry and Alain Rapaport for illuminating discussions.
Malisoff was supported by NSF/DMS Grant 0424011.  De Leenheer was
supported by NSF/DMS Grant 0500861.  The authors thank the referees
for their comments, and they thank Hairui Tu for helping with the
graphics. The second author thanks Ilana Aldor for stimulating
discussions.

\medskip\medskip

 \label{end}
\end{document}